\documentclass{amsart}

\usepackage{amssymb}

\begin{document}

\title[Corrigenda Am. Math. Monthly 92 (1985) 449]{Corrigenda to ``Interesting Series Involving the Central Binomial Coefficient'' [Am. Math. Monthly vol 92 (1985)]}

\author{Richard J. Mathar}
\urladdr{http://www.strw.leidenuniv.nl/~mathar}
\email{mathar@strw.leidenuniv.nl}
\address{Leiden Observatory, Leiden University, P.O. Box 9513, 2300 RA Leiden, The Netherlands}


\date{\today}

\begin{abstract}
These are seven corrigenda to equations in the Lehmer article in
American Mathematical Monthly 92 (1985), pp 449--457, partially
reproduced in the Apelblat tables of integrals and series.
\end{abstract}

\maketitle

\section{Paper by D. H. Lehmer}
The corrigenda to equations in an article by the late
Derrick Lehmer \cite{LehmerAMM92} are:
\begin{enumerate}
\item
Multiply the right hand side of
the first equation on page 451 by two:
\[
\sum_{n=1}^\infty \frac{(-1)^{n+1}{2n\choose n}}{n4^n}= 2\log\left[\frac{\sqrt{2}+1}{2}\right].
\]
The numerical value is sequence A157699 in the Online Encyclopedia of Integer Sequences
(OEIS) \cite{EIS}.

\item
Switch the sign of the first but last term in equation (7):
\[
x\sum_{n=1}^\infty
\frac{{ 2n\choose n}}{n(n+1)}x^n = 2x \log\frac{1-\sqrt{1-4x}}{x}+\frac{\sqrt{1-4x}}{2}
-x(\log 4 -1)-\frac{1}{2}.
\]
\item
Multiply the right hand side of
an equation on page 453 by two:
\[
\sum_{m=1}^\infty \frac{(-1)^{m-1}}{m{2m\choose m}}= 2\log\left(\frac{1+\sqrt{5}}{2}\right)/\sqrt{5}.
\]
The numerical value is sequence A086466 in the OEIS \cite{EIS}.
\item
Flip a sign on one side of another equation  on page 453:
\[
\sum_{m=1}^\infty 
\frac{(-1)^{m-1}}{m^2{2m\choose m}}= 2\left\{ \log\left(\frac{\sqrt{5}+1}{2}\right)\right\}^2.
\]
The numerical value is sequence A086467 in the OEIS \cite{EIS}.
\item
Multiply the right hand side of the last line of another formula on page 453 by four:
\[
\sum_{m=1}^\infty 
\frac{1}{m^3{2m\choose m}}= -\frac{4\zeta(3)}{3}-\frac{\pi\sqrt{3}}{18}\left\{
\psi\left(\frac{1}{3}\right)-\psi\left(\frac{2}{3}\right)
\right\}.
\]
Variants in notation have appeared since then, see theorem 3.3 of a paper
by Borwein et al.\ \cite{BorweinEM10}, which 
is
equation (47) in a later work \cite{BorweinAM70},
or equation (2.67) by Davydychev and Kalmykov
\cite{DavydychevNPB699}.
The numerical value is sequence A145438 in the OEIS \cite{EIS}.
\item
Drop a factor two on the right hand side of a formula on page 455
and change sign on either side:
\[
\sum\frac{(-1)^m m^3}{{2m \choose m}}
=
\frac{2}{625}[14\sigma+5]
.
\]
The numerical value is sequence A157701 \cite{EIS}.
\item
The last equation on page 457 cannot be reproduced;
the value of the sum is $0.280851790115\ldots$ whereas $\pi^2/36=\zeta(2)/6\approx 0.27415567\ldots$
Equation (13) of page 453 yields the substitute
\[
\sum\frac{2^m(2-\sqrt{3})^m}{m^2{2m\choose m}}
=
2\left(\arcsin \tau\right)^2,\quad
\tau \equiv \frac{\sqrt{3}-1}{2} = \sqrt{2}\sin\frac{\pi}{12} =
\sin\frac{\pi}{3}-\sin\frac{\pi}{6}
.
\]
\end{enumerate}

\section{Apelblat Tables}

Equivalent corrections concern
section 4.1 of the Apelblat tables \cite{Apelblat2}:
\begin{enumerate}
\item
Insert an exponent 2 on the left hand side of formula 40:
\[
\sum_{n=1}^\infty \frac{(-1)^{n-1}n^2(n!)^2}{(2n)!} =
\frac{4}{125}\left[ 5-\sqrt{5}\ln\left(\frac{\sqrt{5}+1}{2}\right)\right]
.
\]
The numerical value is sequence A145434 \cite{EIS}.
\item
Change two numbers on the right hand side of formula 42:
\[
\sum_{n=1}^\infty \frac{(-1)^{n-1}(n!)^2}{n(2n)!} =
\frac{2}{\sqrt{5}}\ln\left(\frac{1+\sqrt{5}}{2}\right)
.
\]
\item
Multiply the right hand side of formula 47 by four and
show that the $\psi$ are trigamma functions:
\[
\sum_{n=1}^\infty \frac{(n!)^2}{n^3(2n)!} =
4\left\{
\frac{\pi\sqrt{3}}{72}\left[\psi^{(1)}\left(\frac{2}{3}\right)-\psi^{(1)}\left(\frac{1}{3}\right)\right]-\frac{1}{3}\zeta(3)
\right\}
.
\]
\item
Replace the right hand side of formula 81:
\[
\sum_{n=1}^\infty \frac{2^n(2-\sqrt{3})^n(n!)^2}{n^2(2n)!}=2\arcsin^2\left(
\frac{\sqrt{3}-1}{2}
\right)
.
\]

\end{enumerate}

\bibliographystyle{amsplain}
\bibliography{all}

\end{document}